\documentclass[12pt,amsymb,fullpage]{amsart}
\usepackage{amssymb,amscd,pstricks}

\newtheorem{theorem}{Theorem}[section]
\newtheorem{defn}[theorem]{Definition}

\newtheorem{lemma}[theorem]{Lemma}
\newtheorem{fact}[theorem]{Fact}

\newtheorem{eple}[theorem]{Example}
\newtheorem{rmk}[theorem]{Remarks}
\newtheorem{dsc}[theorem]{Discussion}
\newtheorem{nota}[theorem]{Notation}

\newsavebox{\indbin}
\savebox{\indbin}{\begin{picture}(0,0)
\newlength{\gnu}
\settowidth{\gnu}{$\smile$} \setlength{\unitlength}{.5\gnu}
\put(-1,-.65){$\smile$} \put(-.25,.1){$|$}
\end{picture}}

\newcommand{\be}{\begin{enumerate}}
\newcommand{\bd}{\begin{defn}}
\newcommand{\bt}{\begin{theorem}}
\newcommand{\bl}{\begin{lemma}}
\newcommand{\ee}{\end{enumerate}}
\newcommand{\ed}{\end{defn}}
\newcommand{\et}{\end{theorem}}
\newcommand{\el}{\end{lemma}}

\begin{document}
\title{A Non-Standard Bezout Theorem}
\author{Tristram de Piro}
\address{Mathematics Department, The University of Edinburgh, Kings Buildings, Mayfield Road, Edinburgh, EH9 3JZ} \email{depiro@maths.ed.ac.uk}
\thanks{The author was supported by the William Gordon Seggie Brown research fellowship}
\begin{abstract}
This paper provides a non-standard analogue of Bezout's theorem.
This is acheived by showing that in all characteristics, the
notion of Zariski multiplicity coincides with intersection
multiplicity when we consider the full families of projective
degree d and degree e curves in $P^{2}(L)$. The result is
particularly interesting in that it holds even when we consider
intersections at singular points of curves or when the curves
contain non-reduced components. The proof also provides motivation
for the fact that tangency is a definable relation for families of
curves inside a non-linear $1$-dimensional Zariski structure $X$.
This is a crucial ingredient in unpublished work \cite{Pez} that
any such Zariski structure interprets a pure algebraically closed
field $L$ with $X$ as a definable finite cover.
\end{abstract}

\maketitle

The techniques of non-standard analysis, originally developed for
the real numbers, were recently introduced by Zilber in the
context of Zariski structures. These methods have become extremely
useful to model theorists in answering the question of Zilber's
trichotomy for a large class of strongly minimal sets. This paper
sets out to show that non-standard analysis can also be useful in
algebraic geometry by providing a link with the extensive
machinery developed by model theorists for analytic structures. We
assume some familiarity with certain notions from algebraic and
analytic geometry, as well as the material from Sections 1-5 of
\cite{deP}.
We summarise the relevant facts for the proof in the following three sections;\\

\begin{section}{Etale Morphisms and Algebraic Multiplicity}
\begin{defn}

A morphism $f$ of finite type between varieties $X$ and $Y$ is said to be etale if for all $x\in X$ there are open affine neighborhoods $U$ of $x$ and $V$ of $f(x)$ with $f(U)\subset V$ such that restricted to these neighborhoods the pull back on functions is given by the inclusion;\\

$f^{*}:L[V]\rightarrow L[V]{[x_{1},\ldots, x_{n}]\over f_{1},\ldots, f_{n}}$\\

and $det({\partial f_{i}\over \partial x_{j}})(x)\neq 0\ ,(*)$\\

\end{defn}

The coordinate free definition of etale is that $f$ should be flat
and unramified, where a morphism $f$ is unramified if the sheaf of
relative differentials $\Omega_{X/Y}=0$, clearly this last
confition is satisfied using the condition $(*)$. If we tensor the
exact sequence,\\

 $f^{*}\Omega_{Y}\rightarrow \Omega_{X}\rightarrow
\Omega_{X/Y}\rightarrow 0$\\

with $L(x)$ the residue field of $x$, we obtain an isomorphism\\

$f^{*}\Omega_{Y}\otimes L(x)\rightarrow \Omega_{X}\otimes L(x)$.\\

Identifying $\Omega_{X}\otimes L(x)$ with $T_{x,X}^{*}$ gives
that\\

 $df:(m_{x}/m_{x}^{2})^{*}\rightarrow
(m_{f(x)}/m_{f(x)}^{x})^{*}$\\

 is an isomorphism of tangent spaces or dually $f^{*}(m_{f(x)})=m_{x}$.
 Call this property of etale morphisms $(**)$.\\

We will also require some facts about the etale topology on an
algebraic variety $Y$. We consider a category $Y_{et}$ whose
objects are etale morphisms $U\rightarrow Y$ and whose arrows are
$Y$-morphisms from $U\rightarrow V$. This category has the
following $2$ desirable properties. First given $y\in Y$, the set
of objects of the form $(U,x)\rightarrow (Y,y)$ form a directed
system, namely $(U,x)\subset (U',x')$ if there exists a morphism
$U\rightarrow U'$ taking $x$ to $x'$. Secondly, we can take
``intersections'' of open sets $U_{i}$ and $U_{j}$ by considering
$U_{ij}=U_{i}\times_{Y}U_{j}$; the projection maps are easily show
to be etale and the composition of etale maps is etale, so
$U_{ij}\rightarrow Y$ still lies in $Y_{et}$. If $Y$ is an
irreducible variety over $L$, then all etale morphisms into $Y$
must come from reduced schemes of finite type over $L$, though
they may well fail to be irreducible considered as algebraic
varieties. Now we can define the local ring of $Y$ in the etale toplogy to be;\\

$O_{y,Y}^{\wedge}=lim_{\rightarrow, y\in U}O_{U}(U)$\\

As any open set $U$ of $Y$ clearly induces an etale morphism $U\rightarrow_{i}Y$ of inclusion, we have that $O_{y,Y}\subset O_{y,Y}^{\wedge}$. We want to prove that $O_{y,Y}^{\wedge}$ is a Henselian ring and in fact the smallest Henselian ring containing $O_{y,Y}$. We need the following lemma about Henselian rings;\\

\begin{lemma}

Let $R$ be a local ring with residue field $k$. Suppose that $R$ satisfies the following condition;\\

If $f_{1},\ldots f_{n}\in R[x_{1},\ldots x_{n}]$ and $\bar f_{1}\ldots \bar f_{n}$ have a common root $\bar a$ in $k^{n}$, for which $Jac(\bar f)(\bar a)=({\partial \bar f_{i}\over \partial x_{j}})_{ij}(\bar a)\neq 0$, then $\bar a$ lifts to a common root in $R^{n}$\ (*).\\

Then $R$ is Henselian.\\

\end{lemma}

It remains to show that $O_{y,Y}^{\wedge}$ satisfies $(*)$.

\begin{proof}

Given $f_{1,}\ldots f_{n}$ satisfying the condition of $(*)$, we
can assume the coefficients of the $f_{i}$ belong to
$O_{U_{i}}(U_{i})$ for covers $U_{i}\rightarrow Y$; taking the
intersection $U_{1\ldots i\ldots n}$ we may even assume the
coefficients define functions on a single etale cover $U$ of $Y$.
By the remarks above we can consider $U$ as an algebraic variety
over $K$, and even an affine algebraic variety after taking the
corresponding inclusion. We then consider the variety $V\subset
U\times A^{n}$ defined by $Spec({R(U)[x_{1},\ldots,x_{n}]\over
f_{1},\ldots f_{n}})$. Letting $u\in U$ denote the point in $U$
lying over $y\in Y$, the residue of the coefficients of the
$f_{i}$ at $u$ corresponds to the residue in the local ring $R$,
which tells us exactly that the point $(u,\bar a)$ lies in $V$. By
the Jacobian condition, we have that the projection
$\pi:V\rightarrow U$ is etale at the point $(u,\bar a)$, and hence
on some open neighborhood of $(u,\bar a)$, using Nakayama's Lemma
applied to $\Omega_{V/U}$. Therefore, replacing $V$ by the open
subset $U'\subset V$ gives an etale cover of $U$ and therefore of
$Y$, lying over $y$. Now clearly the coordinate functions
$x_{1},\ldots x_{n}$ restricted to $U'$ lie in
  $O_{y,Y}^{\wedge}$ and lift the root $\bar a$ to a root in $O_{y,Y}^{\wedge}$\\
\end{proof}

We define the Henselization of a local ring $R$ to be the smallest Henselian ring $R'\supset R$, with $R'\subset Frac(R)^{alg}$. We have in fact that;\\

\begin{theorem}

Given an algebraic variety $Y$, $O_{y,Y}^{\wedge}$ is the
Henselization of $O_{y,Y}$

\end{theorem}

The following theorem requires some knowledge of Zariski
structures, see \cite{deP} sections 1-4, or section 2 of this
paper.

\begin{theorem}{Zariski multiplcity is preserved by etale morphisms}

Let $\pi:X\rightarrow Y$ be an etale morphism with $Y$ smooth,
then any $(ab)\in graph(\pi)\subset X\times Y$ is unramified in
the sense of Zariski structures.

\end{theorem}

For this we need the following fact whose algebraic proof relies on the fact that etale morphisms are flat, see \cite{Milne};\\

\begin{fact}

Any etale morphism can be locally presented  in the form \\

\begin{eqnarray*}
\begin{CD}
V@>g>>Spec((A[T]/f(T))_{d})\\
@VV\pi V  @VV\pi' V\\
U@>h>>Spec(A)\\
\end{CD}
\end{eqnarray*}

where $f(T)$ is a monic polynomial in $A[T]$, $f'(T)$ is invertible in $(A[T]/f(T))_{d}$ and $g,h$ are isomorphisms. \\
\end{fact}

Using Lemma 4.6 of \cite{deP} and the fact that the open set $V$ is smooth, we may safely replace $graph(\pi)$ by $\overline {graph (\pi')}\subset F''\times F$ where $F''$ is the projective closure of $Spec((A[T]/f(T))$, $F$ is the projective closure of $Spec(A)$ and $\overline {graph(\pi')}$ is the projective closure of $graph(\pi')$ and show that $(g(b)a)$ is Zariski unramified. Note that over the open subset $U=Spec(A)\subset F$, $\overline{graph(\pi')}=Spec(A[T]/f(T))$ as this is closed in $U\times F''$.  For ease of notation, we replace $(g(b)a)$ by $(ba)$.\\

Suppose that $f$ has degree $n$. Let $\sigma_{1}\ldots \sigma_{n}$ be the elementary symmetric functions in $n$ variables $T_{1},\ldots T_{n}$. Consider the equations\\

$\sigma_{1}(T_{1},\ldots, T_{n})=a_{1}$\\

$\ldots$\\

$\sigma_{n}(T_{1},\ldots,T_{n})=a_{n}$ (*)\\

where $a_{1},\ldots a_{n}$ are the coefficients of $f$ with
appropriate sign. These cut out a closed subscheme $C\subset
Spec(A[T_{1}\ldots T_{N}])$. Suppose $(ba)\in
graph(\pi')=Spec(A[T]/f(T))$ is ramified in the sense of Zariski
structures, then I can find $(a'b_{1}b_{2})\in {\mathcal V}_{abb}$
with $(a'b_{1})$,$(a'b_{2})\in Spec(A(T)/f(T))$ and $b_{1},b_{2}$
distinct. Then complete $(b_{1}b_{2})$ to an $n$-tuple
$(b_{1}b_{2}c_{1}'\ldots c_{n-2}')$ corresponding to the roots of
$f$ over $a'$. The tuple $(a'b_{1}b_{2}c_{1}'\ldots c_{n-2}')$
satisfies $C$, hence so does the specialisation $(abbc_{1}\ldots
c_{n-2})$. Then the tuple $(bbc_{1}\ldots c_{n-2})$ satisfies
$(*)$ with the coefficients evaluated at $a$. However such a
solution is unique up to permutation and corresponds to the roots
of $f$ over $a$. This shows that $f$ has a double root at $(ab)$
and therefore $f'(T)|_{ab}=0$. As $(ab)$ lies inside
$Spec(A[T]/f(T))_{d}$, this contradicts the fact that $f'$ is
invertible in $A[T]/f(T))_{d}$.\\

We also review some facts about algebraic multiplicity and show that algebraic multiplicity
is preserved by etale morphisms.\\

\begin{defn}

Given projective varieties $X_{1}$, $X_{2}$ and a finite morphism
$f:X_{1}\rightarrow X_{2}$, the algebraic multiplicity
$mult_{af(a)}^{alg}(X_{1}/X_{2})$ of $f$ at $a\in X_{1}$ is
$length(O_{a,X_{1}}/f^{*}m_{f(a)})$ where $m_{f(a)}$ is the
maximal ideal of the local ring $O_{f(a)}$.

\end{defn}

\begin{rmk}

Note that this is finite, by the fact that finite morphisms have
finite fibres and the ring $O_{a,X_{1}}/f^{*}m_{f(a)}$ is a
localisation of the fibre $f^{-1}(f(a))\cong
R(f^{-1}(U))\otimes_{R(U)}L\cong R(f^{-1}(U))/m_{f(a)}$ where $U$
is an affine subset of $X_{2}$ containing $f(a)$.

\end{rmk}

We now have the following;\\

\begin{theorem}{Algebraic multiplicity is preserved by etale morphisms};\\

Given finite morphisms $f:X_{3}\rightarrow X_{2}$ and
$g:X_{2}\rightarrow X_{1}$ with $f$ etale. If $a\in X_{3}$, then
$mult_{a,gf(a)}^{alg}(X_{3}/X_{1})=mult_{f(a),gf(a)}^{alg}(X_{2}/X_{1})$.

\end{theorem}

\begin{proof}

This result is essentially given in \cite{Mum}. Let
$O_{f(a),X_{2}}^{\wedge}$ be the Henselisation of the local ring
at $f(a)$. By base change, we have an etale morphism
$f':X'=X_{3}\times_{X_{2}}Spec(O_{f(a)}^{\wedge}) \rightarrow
Spec(O_{f(a)}^{\wedge})$. By the definition of an etale morphism
given above, we may write this cover locally in the form
$Spec(O_{f(a)}^{\wedge}{[x_{1},\ldots,x_{n}]\over f_{1},\ldots,
f_{n}})$, with $det({\partial f_{i}\over \partial x_{j}})\neq 0$
at each closed point in the fibre over $f(a)$. At the closed point
$a$, let $a_{i}$ be the residues of the $x_{i}$ in $L$, then we
have that $(a_{1},\ldots a_{n})$ is a common root for $\{\bar
f_{1},\ldots,\bar f_{n}\}$ where $\bar f_{i}$ is obtained by
reducing $f_{i}$ with respect to the maximal ideal
$m_{f(a),X_{2}}$ of $O_{f(a)}^{\wedge}$. As $O_{f(a)}^{\wedge}$ is
Henselian, by the above, and the determinant condition, we can
lift the roots $a_{i}$ to roots $\alpha_{i}$ of the $f_{i}$ in
$O_{f(a)}^{\wedge}$. We therefore obtain a subscheme
$Z=Spec(O_{f(a)}^{\wedge}{[x_{1},\ldots,x_{n}]\over
<x_{1}-\alpha_{1},\ldots, x_{n}-\alpha_{n}}>)$ of $X'$ which is
isomorphic to $Spec(O_{f(a)}^{\wedge})$ under the restriction of
$f$. Let $Q$ be the $O_{X'}$ ideal defining $Z$, we then have that
$m_{a,X'}=f^{*}m_{f(a),X_{2}}\oplus Q_{a}$. As $f$ is etale, by
$(**)$ after definition 1.1 above, $m_{a,X'}=f^{*}m_{f(a),X_{2}}$,
therefore $Q_{a}=0$ and by Nakayama's lemma $Q=0$ in an open
neighborhood of $a$ in $X'$. This gives that $Z=X'$ in an open
neighborhood of $a$. Hence we obtain the sequence
$O_{f(a),X_{2}}\rightarrow_{f^{*}} O_{a,X_{3}}\rightarrow_{i^{*}}
O_{a,X'}$ (***) where the map $i^{*}f^{*}$ is the inclusion of
$O_{f(a),X_{2}}$ inside $O_{f(a),X_{2}}^{\wedge}$. Now if
$n\subset m_{f(a),X_{2}}$ is the pullback $g^{*}m_{gf(a),X_{1}}$,
we have that
$length(O_{f(a),X_{2}}/n)=length(O_{f(a),X_{2}}^{\wedge}/n)$,
hence the result follows by $(***)$ as required.

\end{proof}

\end{section}

\begin{section}{Zariski Multiplicity}

We work in the context of Theorem 3.3 in \cite{deP}. Namely, $W$
(we used the notation $V$ in \cite{deP}) will denote a smooth
projective variety defined over an algebraically closed field $L$,
considered as a Zariski structure with closed sets given by
algebraic subvarieties defined over $L$. All notions connected to
the definition of Zariski multiplicity will come from a fixed
specialisation map $\pi:W(K_{\omega})\rightarrow W(L)$ where
$K_{\omega}$ denotes a "universal" algebraically closed field
containing $L=K_{0}$. We consider $D$ a smooth subvariety of some
cartesian power $W^{m}$ and a finite cover, with respect to
projection onto the first coordinate, $F\subset D\times W^{k}$,
all defined over $L$ (*). This allows us to make sense of Zariski
multiplicity. In general, we can move freely between Zariski
structure notation and algebraic geometry notation. Clearly $(*)$
makes sense algebraically. Conversely, if $X$ and $Y$ denote fixed
projective varieties defined over $L$ with $Y$ smooth and a finite
morphism $f:X\rightarrow Y$ over $L$ is given , then we can reduce
to the situation of $(*)$ by taking $F$ to be $graph(f)\subset
X\times Y$ with the projection map onto the second factor and $W$
to be the corresponding projective space $P^{n}(L)$ where
$X,Y\subset P^{n}(L)$. We can even take $W$ to be the
$1$-dimensional Zariski structure $P^{1}(L)$ by using the
embedding of $P^{n}(L)$ into the
$N$'th Cartesian power of $P^{1}(L)$ for sufficiently large $N$.\\

We use the definition of Zariski multiplicity for irreducible
finite covers given in 4.1 of \cite{deP}. We will also require the
following generalisation;\\

\begin{defn}

Let  $F\subset D\times W^{k}$ be an equidimensional, finite cover
of smooth $D$, with irreducible components $C_{1},\ldots,C_{n}$.
Then for $(ab)\in F$, we define $Mult_{ab}(F/D)=\sum_{(ab)\in
C_{i}}Mult_{ab}(C_{i}/D)$.

\end{defn}

Clearly this is well defined using the definition of Zariski
multiplicity for irreducible covers. However, until Lemma 2.9, the
assumption that $F$ is irreducible will be in force.

\begin{lemma}{Zariski multiplicity is multiplicative over composition}

Suppose that $F_{1},F_{2}$ and $F_{3}$ are smooth, irreducible,
with $F_{2}\subset F_{1}\times W^{k}$ and $F_{3}\subset
F_{2}\times W^{l}$ finite covers. Let $(abc)\in F_{3}\subset
F_{1}\times W^{k}\times W^{l}$. Then
$mult_{abc}(F_{3}/F_{1})=mult_{ab}(F_{2}/F_{1})mult_{abc}(F_{3}/F_{2})$.

\end{lemma}

\begin{proof}
To see this, let $m=mult_{ab}(F_{2}/F_{1})$ and $n=mult_{abc}(F_{3}/F_{2})$. Choose $a'\in {\mathcal V}_{a}\cap F_{1}(K_{\omega})$ generic over $L$. By definition, we can find distinct $b_{1}\ldots b_{m}$ in $W^{k}(K_{\omega})\cap{\mathcal V_{b}}$ such that $F_{2}(a',b_{i})$ holds. As $F_{2}$ is a finite cover of $F_{1}$, we have that $dim(a'b_{i}/L)=dim(a'/L)=dim(F_{1})=dim(F_{2})$, so each $(a'b_{i})\in {\mathcal V}_{ab}\cap F_{2}$ is generic over $L$. Again by definition, we can find distinct $c_{i1}\ldots c_{in}$ in $W^{l}(K_{\omega})\cap{\mathcal V_{c}}$ such that $F_{3}(a'b_{i}c_{ij})$ holds. Then the $mn$ distinct elements $(a'b_{i}c_{ij})$ are in ${\mathcal V_{abc}}$, so by definition of multiplicity $mult_{abc}(F_{3}/F_{1})=mn$ as required.\\
\end{proof}

\begin{lemma}

Let hypotheses be as in the above lemma with the extra condition
that the cover $F_{3}/F_{2}$ is etale. Then for $(abc)\in F_{3}$,
$mult_{abc}(F_{3}/F_{1})=mult_{ab}(F_{2}/F_{1})$

\end{lemma}

\begin{proof}

This is an immediate consequence of Lemma 2.2 and Theorem 1.4.

\end{proof}

\begin{lemma}{Zariski multiplicity is summable over
specialisation}

Suppose that $F\subset D\times W^{k}$ is a finite irreducible
cover with $D$ smooth. Suppose $(ab)\in F$, $a'\in {\mathcal
V_{a}}\cap D$ and $a''\in {\mathcal V_{a'}}\cap D$ with $a''$
generic over $L$.
Then\\

$Mult_{ab}(F/D)=\Sigma_{b'\in {\mathcal V}_{b}\cap
F(a'y)}Mult_{a'b'}(F/D)$\\

\end{lemma}

\begin{proof}

Suppose $F(a''b_{1}),\ldots F(a''b_{n})$ hold with $b_{i}\in
{\mathcal V}_{b}$, so $\{b_{1},\ldots, b_{n}\}$ witness the fact
that $Mult_{ab}(F/D)=n$. Write $\{b_{1},\ldots b_{n}\}$ as
$\{b_{11},\ldots,\\
b_{1m_{1}},b_{21},\ldots, b_{2m_{2}},\ldots,b_{i1},\ldots
b_{ij},\ldots,b_{im_{i}},\ldots, b_{nm_{n}}\}$ (*), where $b_{ij}$
maps to $a_{i}$ in the specialisation taking $a''$ to $a'$. To
prove the lemma, it is sufficient to show that $F(a'y)\cap
{\mathcal V}_{b}=\{a_{1},\ldots,a_{n}\}$ and
$Mult_{(a'a_{i})}(F/D)=m_{i}$. The second statement just follows
from the fact that $a''$ is generic in $D$ over $L$ in ${\mathcal
V}_{a'}$. To prove the first statement, suppose we can find
$a_{n+1}$ with $F(a'a_{n+1})$ and $a_{n+1}\in {\mathcal V}_{b}$
but $a_{n+1}\notin \{a_{1},\ldots a_{n}\}$. By Theorem 3.3 in
\cite{deP}, we can find $c$ with $F(a''c)$ and $(a''c)$
specialising to $(a'a_{n+1})$. As $a_{n+1}\in{\mathcal V}_{b}$,
$(a'a_{n+1})$ specialises to $(ab)$, hence so does $(a''c)$.
Therefore, $c$ must witness the fact that $Mult_{ab}(F/D)=n$ and
appear in the set $\{b_{1},\ldots,b_{n}\}$. This clearly
contradicts the arrangement of $\{b_{1},\ldots,b_{n}\}$ given in
$(*)$.

\end{proof}

\begin{defn}

Let $F\subset U\times V\times W^{k}$ be an irreducible finite
cover of $U\times V$ with $U$ and $V$ smooth.\\

Given $(u,v,x)\in F$ we define;\\

$Left Mult_{u,v,x}(F/D)=Card({\mathcal V}_{x}\cap F(u',v))$ for
$u'\in{\mathcal V}_{u}\cap U$ generic over $L$. \\

$Right Mult_{u,v,x}(F/D)=Card({\mathcal V}_{x}\cap F(u,v'))$ for
$v'\in{\mathcal V}_{v}\cap V$ generic over $L$. \\

\end{defn}

We first show that both left and right multiplicity are well
defined. In order to see this, observe that the fibres $F(u,V)$
and $F(U,v)$ are finite covers of $V$ and $U$ respectively with
$U$ and $V$ smooth. Moreover, the fibres $F(u,V)$ and $F(U,v)$ are
equidimensional covers of $V$ and $U$ respectively. In order to
see this, as $U$ is smooth, it satisfies the presmoothness axiom
with the smooth projective variety $W^{k}$ given in Definition 1.1
of \cite{deP}. The fibre $F(u,V)=F\cap (W^{k}\times \{u\}\times
V)$. By presmoothness, each irreducible component of the
intersection has dimension at least $dim(F)+dim(W^{k}\times
V)-dim(U\times V\times W^{k})=dim(F)-dim(U)=dim(V)$. As $F(u,V)$
is a finite cover of $V$, it has exactly this dimension. Now we
can use the definition
of Zariski multiplicity given in 1.4.\\

We then claim the following;\\

\begin{lemma}{Factoring Multiplicity}

In the situation of the above definition, we have that;\\

$Mult_{u,v,x}(F/U\times V)=\Sigma_{x'\in ({\mathcal V}_{x}\cap
F(y,u',v))} Right Mult_{x',u',v}(F/U\times V)$ for $u'$ generic in $U$ over $L$.\\

$Mult_{u,v,x}(F/U\times V)=\Sigma_{x'\in ({\mathcal V}_{x}\cap
F(y,u,v'))} Left Mult_{x',u,v'}(F/U\times V)$ for $v'$ generic in $V$ over $L$.\\

\end{lemma}

\begin{proof}

We just prove the first statement, the proof of the second is
apart from notation identical. By the construction in section 2
and Lemma 3.2 of \cite{deP}, we can choose algebraically closed
fields $L=K_{0}\subset K_{n_{1}}\subset K_{n_{2}}\subset
K_{\omega}$, and tuples $u'\in K_{n_{1}}$, $v'\in K_{n_{2}}$ such
that $u'$ is generic in $U$ over $L$, $v'$ is generic in $V$ over
$K_{n_{1}}$ with specialisations $\pi_{1}:
P^{n}(K_{n_{1}})\rightarrow P^{n}(L)$ and
$\pi_{2}:P^{n}(K_{n_{2}})\rightarrow P^{n}(K_{1})$ such that
$\pi_{2}(u'v')=(u'v)$ and $\pi_{1}(u'v)=(uv)$. Now
$dim(u'v'/L)=dim(v'/L(u'))+dim(u'/L)=dim(V)+dim(U)$, hence $u'v'$
is generic in $U\times V$ over $L$. Therefore
$Mult_{u,v,x}=Card({\mathcal V}_{x}\cap F(u'v'))$. Let
$S=\{y_{11},\ldots,y_{1m_{1}},\ldots,y_{ij_{i}},\ldots,
y_{n1},\\
\ldots,y_{nm_{n}}\}$ be distinct elements in ${\mathcal V}_{x}\cap
W^{k}$ witnessing this multiplicity such that for $1\leq j_{i}\leq
m_{i}$, $\pi_{2}(y_{ij_{i}})=z_{i}\in {\mathcal V_{x}}\cap W^{k}$.
It is sufficient to show that $Right Mult_{u'v,z_{i}}(F/U\times
V)=m_{i}$ and $\{z_{1},\ldots z_{n}\}$ enumerates ${\mathcal
V}_{x}\cap F(y,u',v)$. The first statement follows as
$v'\in{\mathcal V}_{v}\cap V$ is generic in $V$ over $L(u')$. For
the second statement, suppose that we can find $z_{n+1}\in
{\mathcal V}_{x}\cap F(y,u',v)$ with $z_{n+1}\notin \{z_{1},\ldots
z_{n}\}$. Consider $F(u',V)$ as a finite cover of $V$, defined
over $L(u')$, so by the above $F(u',V)$ is an equidimensional
finite cover of $V$. Then, as $v'$ was chosen to be generic in $V$
over $L(u')$, choosing an irreducible component of $F(u',V)$
passing through $(z_{n+1},u'v)$, by the lifting result of Theorem
3.3 in \cite{deP}, we can find $y_{n+1}\in {\mathcal
V}_{z_{n+1}}\cap W^{k}$ such that $F(y_{n+1},u',v')$. Clearly,
$y_{n+1}\in S$ which contradicts the definition of $S$.

\end{proof}

Theorem 3.3 of \cite{deP} does not hold in the case when $D$ fails
to be smooth. However, in the case of etale covers, we still have
the following result;

\begin{lemma}{Lifting Lemma for Etale Covers}

Let $F\subset D\times W^{k}$ be an etale cover of $D$ defined over
$L$, with the projection map denoted by $f$. Then given $a \in D$,
$(ab)\in F$ and $a'\in {\mathcal V}_{a}\cap D$ generic over $L$,
we can find $b'\in {\mathcal V}_{b}$ such that $F(a',b')$ holds.
Moreover $b'$ is unique, hence $Mult_{ab}(F/D)=1$. Moreover, in
the situation of Lemma 2.3, without requiring that $F_{2}$ is
smooth, we have that for $(abc)\in F_{3}$,
$mult_{abc}(F_{3}/F_{1})=mult_{ab}(F_{2}/F_{1})$.

\end{lemma}

\begin{proof}

Using the definition of etale given in section 1 above, we can
assume that the cover is given algebraically in the form
$f^{*}:L[D]\rightarrow L[D]{[x_{1},\ldots,x_{n}]\over
f_{1},\ldots,f_{n}}$ with $det({\partial f_{i}\over \partial
x_{j}})_{ij}(x)\neq 0$ for all $x\in F$. So we can present the
cover in the form $f_{1}(x,y)=0, f_{2}(x,y)=0,\ldots,
f_{n}(x,y)=0$, with $y$ in $D$ and $x$ in $A^{n}(L)$. Let $L_{m}$
be the algebraic closure of the field generated by $L$ and $\bar g
(a)$ where $\bar g$ is a tuple of functions defining $D$ locally.
Consider the system of equations $f_{1}(x,a)=f_{2}(x,a)=\ldots
=f_{n}(x,a)=0$ defined over $L_{m}$. Then this system is solved by
$b$ in $L_{m}$ with the property that $det({\partial
f_{i}\over\partial x_{j}})_{ij}(b)\neq 0$ (*). Now suppose that
$a'\in {\mathcal V_{a}}\cap D$ is chosen to be generic over $L$.
By the construction given in $2$ of \cite{deP} and the following
Lemma 2.2, we may assume that $a'$ lies in $L_{s}[[t^{1/r}]]$, the
formal power series in the variable $t^{1/r}$ for some
algebraically closed field $L_{s}$ extending $L_{m}$. This is a
henselian ring, hence if we consider the system of equations
$f_{1}(x,a')=f_{2}(x,a')=\ldots=f_{n}(x,a')=0$ with coefficients
in $L_{s}[[t^{1/r}]]$, by the fact that the system specialises to
a solution in $L_{s}$ with the condition (*) we can find a
solution $b'$ in $L_{s}[[t^{1/r}]]$. Then $(a'b')$ lies in $F$ and
by construction $b'\in{\mathcal V}_{b}$. The uniqueness result
follows from the proof of Theorem 1.4. For the last part, suppose
that $mult_{ab}(F_{2}/F_{1})=n$, then we can find $a'\in{\mathcal
V_{a}}\cap F_{1}$ generic over $L$ and $\{b_{1},\ldots b_{n}\}\in
{\mathcal V_{b}}\cap W^{k}$ distinct such that $F(a',b_{i})$
holds. Each $(a'b_{i})$ is generic in $F_{2}$ over $L$, hence by
the previous part of the lemma, we can find a unique $c_{i}\in
{\mathcal V_{c}}\cap W^{l}$ such that $F_{3}(a'b_{i}c_{i})$ holds.
This show that $mult_{abc}(F_{3}/F_{1})=n$ as required.

\end{proof}

\begin{lemma}{Lifting Lemma for Etale Covers with Right(Left) Multiplicity}

Let hypotheses be as in Lemma 1.5, with the additional assumption
that $F_{1}=U\times V$, $F_{2}$ is a smooth irreducible cover of
$F_{1}$ and $F_{3}$ is an irreducible etale cover of $F_{2}$. Then
with, notion as in the lemma, given $(uvbc)\in F_{3}$, $Right
Mult_{uvbc}(F_{3}/F_{1})=Right Mult_{uvb}(F_{2}/F_{1})$. Similarly
for left multiplicity.

\end{lemma}

\begin{proof}

Suppose that $Right Mult_{uvb}(F_{2}/F_{1})=n$, then for
$v'\in{\mathcal V_{b}}$ generic in $V$ over $L$, we can find
$\{b_{1},\ldots,b_{i},\ldots b_{n}\}\in{\mathcal V_{b}}$ with
$F_{2}(uv'b_{i})$ holding. For each $b_{i}$ we claim that there
exists a unique $c_{i}\in {\mathcal V}_{c}$ such that
$F_{3}(uv'b_{i}c_{i})$ holds. For the existence, we can use Lemma
2.7, with the simple modification that, with the notation there,
if $L_{m}$ is the algebraic closure of the field generated by
${\bar g}(uv)$, then provided $dim(V)\geq 1$, we can find
$v'\in{\mathcal V}_{v}\cap V$ generic over $L$ with $uv'\in
L_{s}[[t^{1/r} ]]$ for some algebraically closed field $L_{s}$
containing $L_{m}$. For the uniqueness, we can use the fact that
Zariski multiplicity is summable over specialisation (Lemma 2.4)
and the fact that for generic $(u'v'b_{i}')\in {\mathcal
V_{uvb}}\cap F_{2}$, we can find a unique $c_{i}'\in{\mathcal
V}_{c}$ such that $F_{3}(u'v'b_{i}'c_{i}')$ holds. Finally, we
claim that $\{b_{1}c_{1},\ldots, b_{n}c_{n}\}$ enumerate
$F_{3}(uv'xy)\cap {\mathcal V}_{bc}$. This is clear by the above
proof and the fact that $\{b_{1},\ldots,b_{n}\}$ enumerates
$F_{2}(uv'x)\cap {\mathcal V}_{b}$.

\end{proof}

\begin{lemma}

The following versions of the above properties hold when we
consider finite equidimensional covers, possibly with components,
with the definition of Zariski multiplicity given in 2.1.

\end{lemma}

\begin{proof}

For Lemma 2.3, we replace the hypotheses with $F_{1}$ is smooth
irreducible, $F_{2}$ is an equidimensional finite cover of $F_{1}$
and $F_{3}$ is an etale cover of $F_{2}$. We then claim, using
notation as in Lemma 2.2, that
$mult_{abc}(F_{3}/F_{1})=mult_{ab}(F_{2}/F_{1}$. By definition
$mult_{abc}(F_{3}/F_{1})=\sum_{(abc)\in
C_{i}}(mult_{abc}(C_{i}/F_{1}))$, where $C_{i}$ are the
irreducible components of $F_{3}$ passing through $(abc)$. As
$F_{3}$ is an etale cover of $F_{2}$, the images of the $C_{i}$
are precisely the irreducible components $D_{i}$ of $F_{2}$
passing through $(ab)$, each $C_{i}$ is an etale cover of $D_{i}$
and $mult_{ab}(F_{2}/F_{1})=\sum_{(ab)\in
D_{i}}(mult_{ab}(D_{i}/F_{1}))$. Hence, it is sufficient to prove
the result in the case when $F_{2}$ and $F_{3}$ are irreducible.
This is just Lemma 2.3\\

For Lemma 2.4, we replace the hypothesis with $F$ is an
equidimensional finite cover of $D$. The proof then goes through
exactly as in the lemma with the observation that if we find
$a_{n+1}\in{\mathcal V}_{b}$ and $F(a'a_{n+1})$ then we can find
an irreducible component $C$ passing through $(a'a_{n+1})$ which
allows us to apply Theorem 3.3 in \cite{deP} to obtain $c$ with
$C(a''c)$ and $(a''c)$ specialising to $(a'a_{n+1})$.\\

For Definition 2.5, we alter the hypothesis to $F$ is an
equidimensional finite cover of $U\times V$. Again, we can use an
identical proof to show that left multiplicity and right
multiplicity are well defined. The proof of Lemma 2.6 with the
new hypothesis on $F$ is identical\\

We don't require a modified version of Lemma 2.7, the result we
need is contained in the modified proof of Lemma 2.3\\

For Lemma 2.8, we alter the hypotheses to $F_{2}$ is an
equidimensional cover of $F_{1}$ and $F_{3}$ is an etale cover of
$F_{2}$. We then claim that for $(uvb)$ a non-singular point of
$F_{2}$ and $(uvbc)\in F_{3}$, necessarily non-singular as well,
that $Right Mult_{uvbc}(F_{3}/F_{1})=Right
Mult_{uvb}(F_{2}/F_{1})$ and similarily for left multiplicity. To
prove this, note that as $(uvb)$ and $(uvbc)$ are non-singular
points, there exist unique components $C$ and $D$ passing through
$(uvb)$ and $(uvbc)$ respectively. Now replacing $C$ and $D$ by
the open subsets $C'$ and $D'$ of smooth points, we can apply the
definition of Right Multiplicity and the proof of Lemma 2.8.

\end{proof}

\end{section}

\begin{section}{Analytic Methods}

In order to use the method of etale morphisms, which preserve
Zariski multiplicity, we need to work inside the Henselisation of
local rings $L[x_{1},\ldots,x_{n}]_{(x_{1},\ldots,x_{n})}$. In the next section, we will
only need the result for the local ring in $2$ variables $L[x,y]_{(x,y)}$.\\

We let $L[[x_{1},\ldots,x_{n}]]$ denote the ring of formal power
series in $n$ variables, which is the formal completion of
$L[x_{1},\ldots,x_{n}]_{(x_{1},\ldots,x_{n})}$ with respect to the
canonical order valuation, see for example Section 2 of
\cite{deP}. The following is a classical result, requiring the
fact that etale morphisms are flat, used in the proof of the Artin
approximation theorem. This relates the Henselisation of the ring
$L\{x_{1},\ldots,x_{n}\}$ of strictly convergent power series in
several variables with its formal completion
$L[[x_{1},\ldots,x_{n}]]$, see \cite{Art} or \cite{Rob};\\

Henselisation$(L[x_{1},\ldots x_{n}]_{(x_{1},\ldots x_{n})})=L[[x_{1},\ldots x_{n}]]\cap L(x_{1},\ldots x_{n})^{alg}$\\

This implies that\\

$O_{\bar 0, A^{n}}^{\wedge}\cong L[[x_{1},\ldots x_{n}]]\cap L(x_{1},\ldots x_{n})^{alg}$\\

The following result, which can be found in \cite{bgr}, is essential for the next section\\

\begin{lemma}{Weierstrass Preparation}

Let $F(x_{1},\ldots x_{n})$ be a polynomial in
$L[x_{1},\ldots,x_{n}]$ which is regular in the variable $x_{n}$.
Then we have $F(x_{1},\ldots,
x_{n})=U(x_{1},\ldots,x_{n})G(x_{1},\ldots,x_{n})$ where
$U(x_{1},\ldots,x_{n})$ is a unit in the local ring
$L[[x_{1},\ldots,x_{n}]]$ and \\
$G(x_{1},\ldots,\ldots x_{n})$ is
a Weierstrass polynomial in $x_{n}$ with coefficients in
$L[[x_{1},\ldots,x_{n-1}]]$

\end{lemma}

We will require the Weierstrass decomposition to hold inside\\
Henselisation$(L[x_{1},\ldots,x_{n}])$, therefore we need to show
that the Weierstass data can be found inside
$L(x_{1},\ldots,x_{n})^{alg}$. This is acheived by the following
lemma;\\

\begin{lemma}{Definability of Weierstrass data}

Let $F(x_{1},\ldots,x_{n})$ be a polynomial with coefficients in
$L$ such that $F$ is regular in $x_{n}$, then if
$F(x_{1},\ldots,x_{n})=U(x_{1},\ldots,x_{n})G(x_{1},\ldots,x_{n})$
is the Weierstrass decomposition of $F$ with
$G(x_{1},\ldots,x_{n})=x_{n}^{m}+a_{1}(x_{1},\ldots,x_{n-1})x_{n}^{m-1}+\ldots
+a_{m}(x_{1},\ldots,x_{n-1})$, and $a_{i}\in
L[[x_{1},\ldots,x_{n-1}]]$, $U(x_{1},\ldots,x_{n})\in
L[[x_{1},\ldots,x_{n}]]$, then $a_{i}(x_{1},\ldots,x_{n-1})\in
L(x_{1},\ldots,x_{n-1})^{alg}$ and $U(x_{1},\ldots,x_{n})\in
L(x_{1},\ldots,x_{n})^{alg}$.

\end{lemma}

\begin{proof}

 This can be proved by rigid analytic methods. Equip $L$ with
a complete non-trivial non-archimedean valuation $v$ and
corresponding norm $||.||_{v}$, this can be done for example by
assuming that $L$ is the completion of an algebraically closed
field with any non-archimidean valuation, see \cite{bgr}. Let
$T_{n-1}(L)$ be the free Tate algebra in the indeterminate
variables $x_{1},\ldots,x_{n-1}$ over $L$, that is the subalgebra
of strictly convergent power series in
$L[[x_{1},\ldots,x_{n-1}]]$. By the proof of Weierstrass
preparation in \cite{bgr}, as $F\in T_{n-1}(L)[x_{n}]$, the
coefficients $a_{i}$ lie in $T_{n-1}(L)$ and
$U(x_{1},\ldots,x_{n})\in T_{n-1}(L)[x_{n}]$. Now choose
$(u_{1},\ldots u_{n-1})\subset L$ transcendental over the
coefficients of $F$ with $max(\{||u_{i}||\})\leq 1$. Then if
$s_{1}(\bar u),\ldots, s_{m}(\bar u)$ denote the roots of $F(\bar
u,x_{n})$ with $||s_{i}(\bar u)||\leq 1$, then both $U(\bar
u,s_{i}(\bar u))$ and $G(\bar u,s_{i}(\bar u))$ define elements of
$L$ and moreover, by a theorem in \cite{Rob}, we have that the
coefficients $a_{i}(\bar u)$ are symmetric functions of the
$s_{i}(\bar u)$. Hence the $a_{i}(\bar u)$ belong to $L(\bar
u)^{alg}$. As $\bar u$ was transcendental, we have that each
$a_{i}\in L[x_{1},\ldots,x_{n-1}]^{alg}$. As $U(x_{1},\ldots
x_{n})=F/G(x_{1},\ldots,x_{n})$, we clearly have that
$U(x_{1},\ldots,x_{n})\in L[x_{1},\ldots,x_{n}]^{alg}$ as well.

\end{proof}

\end{section}

\begin{section}{Families of Curves in $P^{2}(L)$}

We consider the family $Q_{d}$ of projective curves in $P^{2}(L)$
with degree $d$. An element of $Q_{d}$ may be
written;\\

$\sum_{0\leq i+j\leq d}a_{ij}(X/Z)^{i}(Y/Z)^{j}=0$\\

which, rewriting in homogenous form, becomes;\\

$\sum_{0\leq i+j\leq d}a_{ij}X^{i}Y^{j}Z^{d-(i+j)}=0$\\

For ease of notation, we will use affine coordinates $x=X/Z$ and
$y=Y/Z$. More generally, if we give an affine cover, we implicitly
assume that it can be projectivized by taking $\bar
y=(y_{1},\ldots,y_{n})=(Y_{1}/Z,\ldots,Y_{n}/Z)$. As the notion of Zariski multiplicity
is local, this will not effect our calculations.\\

Now consider two such families $Q_{d}$ and $Q_{e}$. Then we have
the cover obtained by intersecting degree $d$ and degree $e$ curves\\

$Spec(L[x,y,u_{ij},v_{ij}]/<s(u_{ij},x,y),
t(v_{ij},x,y)>)\rightarrow Spec(L[u_{ij},v_{ij}])$. (*)\\

where\\

$s(u_{ij},x,y)=\sum_{0\leq i+j\leq d}u_{ij}x^{i}y^{j}$\\

$t(v_{ij},x,y)=\sum_{0\leq i+j\leq e}v_{ij}x^{i}y^{j}$\\

We denote the parameter space for degree $d$ curves by $U$ and the
parameter space for degree $e$ curves by $V$. These are affine
spaces of dimension $(d+1)(d+2)/2$ and $(e+1)(e+2)/2$
respectively. The cover (*) is generically finite, that is there
exists an open subset $U'\subset Sp(L[u_{ij},v_{ij}])$ for which
the restricted cover has finite fibres. Throughout this section,
we will denote the base space of the cover by $U\times V$, bearing
in mind that we implicitly mean by this $(U\times V)\cap U'$. Now,
given $2$ fixed parameters sets $\bar u$ and $\bar v$, with $(\bar
u,\bar v)\in U'$, corresponding to curves $C_{\bar u}$ and
$C_{\bar v}$, the algebraic multiplicity of the cover $(*)$ at
$(00,\bar u,\bar v)$ is exactly the intersection multiplicity
$I(C_{\bar u},C_{\bar v},00)$ of the curves at $(00)$. The cover
(*) is equidimensional as $U\times V$ satisfies the presmoothness
axiom with the smooth projective variety $P^{2}(L)$. Restricting
to a finite cover over $U'$, by definition 2.1 we can also define
the Zariski multiplicity of the cover at the point $(00,\bar
u,\bar v)$. The main result that we shall prove in this paper is
the following, which
 generalises an observation given in \cite{Mark};\\

\begin{theorem}

In all characteristics, the algebraic multiplicity and Zariski
multiplicity of the cover $(*)$ coincide at $(00,\bar u,\bar v)$.

\end{theorem}

\begin{defn}

We say that a monic polynomial $p(x,\bar y)$ is Weierstrass in $x$
if $p(x,\bar y)=x^{n}+\ldots +q_{j}(\bar y)x^{n-j}+\ldots
+q_{n}(\bar y)$ with $q_{j}(\bar 0)=0$.

\end{defn}

\begin{defn}

Let $F(x,\bar y)$ be a polynomial in $x$ with coefficients in
$L[\bar y]$. We say the cover\\

$Spec(L[x\bar y]/<F(x,\bar y)>)\rightarrow Spec(L[\bar y])$\\

is generically reduced if for generic $\bar u\in Spec(L[\bar y])$,
$F(x,\bar u)$ has no repeated roots.

\end{defn}

\begin{defn}

Let $F\rightarrow U\times V$ be a finite cover with $U$ and $V$
smooth, such that for $(\bar u,\bar v)\in U\times V$ the fibre
$F(\bar u,\bar v)$ consists of the intersection of algebraic
curves $F_{\bar u},F_{\bar v}$. We call the family sufficiently
deformable at $(\bar u_{0},\bar v_{0})$ if there exists $\bar
u'\in U$ generic over $L$ such that $F_{\bar u'}$ intersects
$F_{\bar v_{0}}$ transversely at simple points.

\end{defn}

We now require a series of lemmas;\\

\begin{lemma}

Let $F(x,\bar y)$ be a Weierstrass polynomial in $x$ with
$F(0,\bar 0)=0$ then algebraic multiplicity and Zariski
multiplicity coincide at $(0,\bar 0)$ if the
cover\\

$Spec(L[x\bar y]/<F(x,\bar y)>)\rightarrow Spec(L[\bar y])$\\

is generically reduced.

\end{lemma}
\begin{proof}
We have that $F(x,\bar y)=x^{n}+q_{1}(\bar
y)x^{n-1}+\ldots+q_{n}(\bar y)$ where $q_{i}(\bar 0)=0$. The
algebraic multiplicity is given by $length(L[x]/F(x,\bar
0))=ord(F(x,\bar 0)=n$ in the ring $L[x]$ with the canonical
  valuation. We first claim that the Zariski multiplicity is the number of solutions to
   $x^{n}+q_{1}(\bar{\epsilon})x^{n-1}+\ldots+q_{n}(\bar{\epsilon})=0$ (\dag), where
    $\bar \epsilon$ is generic in ${\mathcal V}_{\bar 0}$. For suppose that
    $(a,\bar{\epsilon})$ is such a solution, then $F(a,\bar{\epsilon})=0$ and
    by specialisation $F(\pi(a),\bar 0)=0$. As $F$ is a Weierstrass polynomial in $x$,
     $\pi(a)=0$, hence $a\in{\mathcal V}_{0}$, giving the claim.
     We have that $Disc(F(x,\bar y))=Res_{\bar y}(F,{\partial F\over \partial x})$
is a regular polynomial in $\bar y$ defined over $L$. By the fact
that the cover is generically reduced, this defines a proper
closed subset of $Spec(L[\bar y])$. Therefore, $Disc(F(x,\bar
y))|\bar{\epsilon}\neq 0$, hence (\dag) has no repeated roots.
This gives the lemma.
\end{proof}

\begin{lemma}

Let $F(x,\bar y)$ be any polynomial with $F(x,\bar 0)\neq 0$ and
$F(0,\bar 0)=0$. Then if the cover $Spec(L[x,\bar y]/<F(x,\bar
y)>)\rightarrow Spec(L[\bar y])$ is generically reduced, the
Zariski multiplicity at $(0,\bar 0)$ equals $ord(F(x,\bar 0))$ in
$L[x]$.

\end{lemma}
\begin{proof}
By the Weierstrass Preparation Theorem, Lemma 3.1, we can write $F(x,\bar y)=U(x,\bar y)G(x,\bar y)$ with $U(x,\bar y),G(x,\bar y)\in L[[x,\bar y]]$, $G(x,\bar y)$ a Weierstrass polynomial in $x$ and $deg(G)=ord(F(x,\bar 0))$, see also the more closely related statement given in \cite{Aby1}. By Lemma 3.2, we may take the new coefficients to lie inside the Henselized ring $L[x,\bar y]_{\bar 0}^{\wedge}$, hence inside some finite etale extension $L[x,\bar y]^{ext}$ of $L[x,\bar y]$ (possibly after localising $L[x,\bar y]$ corresponding to an open subset of $Spec(L[x,\bar y])$ containing $(0,\bar 0)$). Now we have the sequence of morphisms;\\

$Sp(L[x,\bar y]^{ext}/UG)\rightarrow Spec(L[x,\bar y]/F)\rightarrow Spec(L[\bar y])$\\

The left hand morphism is etale at $\bar 0$, hence by Lemma 2.3 or Lemma 2.7, to compute the Zariski multiplicity of the right hand morphism, we need to compute the Zariski multiplicity of the cover\\

$Spec(L[x,\bar y]^{ext}/UG)\rightarrow Spec(L[\bar y])$\\

at $(0,\bar 0)^{lift}$, the marked point in the cover above
$(0,\bar 0)$. Choose $\bar \epsilon\in {\mathcal V}_{\bar 0}$, the
fibre of the cover is given formally analytically by $L[[x,\bar
y]]/<UG>\otimes_{L[\bar y],\bar y\mapsto\bar\epsilon}L$, hence by
solutions to $U(x,\bar \epsilon)G(x,\epsilon)$. By definition of
Zariski multiplicity, we consider only solutions $(x\bar\epsilon)$
in ${\mathcal V}_{(0,\bar 0)^{lift}}$.
As $U(x,\bar y)$ is a unit in the local ring $L[x,\bar y]^{ext}_{(0,\bar 0)^{lift}}$, we must have $U(x,\bar\epsilon)\neq 0$ for such solutions, otherwise by specialisation $U((0,\bar 0)^{lift})=0$. Hence, the solutions are given by $G(x,\bar \epsilon)=0$. Now, we use the previous lemma to give that the Zariski multiplicity is exactly $deg(G)$ as required.\\

\end{proof}

Now return to the cover\\

$Sp(L[x,y,u_{ij},v_{ij}]/<s(u_{ij},x,y),t(v_{ij},x,y)>)\rightarrow
Sp(L[u_{ij},v_{ij}])$ (*)\\

We will show below, Lemma 4.12, that this is a sufficiently
deformable family at $(\bar u_{0},\bar v_{0})$ when $C_{\bar
u_{0}}$ and $C_{\bar v_{0}}$ define reduced curves.
We claim the following;\\

\begin{lemma}

Suppose parameters $\bar u^{0}$ and $\bar v^{0}$ are chosen such
that $C_{\bar u^{0}}$ and $C_{\bar v^{0}}$ are reduced Weierstrass
polynomials in $x$. Then the Zariski multiplicity of the cover
$(*)$ at $(0,0,\bar u^{0},\bar v^{0})$ equals the intersection
multiplicity $I(C_{\bar u^{0}}, C_{\bar v^{0}},(0,0))$ of $C_{\bar
u^{0}}$ and $C_{\bar v^{0}}$ at $(0,0)$.

\end{lemma}

\begin{proof}

Introduce new parameters $\bar u'$ and $\bar v'$. Let $C_{\bar
u^{0}}^{\bar u'}$ and $C_{\bar v^{0}}^{\bar v'}$ denote the curves
$C_{\bar u^{0}}$ and $C_{\bar v^{0}}$ deformed by the parameters
$\bar u'$ and $\bar v'$ respectively. That is $C_{\bar
u^{0}}^{\bar u'}$ is given by the new equation $\Sigma_{1\leq
i+j\leq d}(u^{0}_{ij}+u_{ij}')x^{i}y^{j}$. Let $F(y,\bar u',\bar
v')=Res(C_{\bar u^{0}}^{\bar u'},C_{\bar v^{0}}^{\bar v'})$.
Then,\\

$F(0,\bar 0,\bar 0)=Res(s(u^{0}_{ij},x,0),t(v^{0}_{ij},x,0))=0$\\

as $C_{\bar u^{0}}$ and $C_{\bar v^{0}}$ are Weierstrass in $x$
and share a common solution at $(0,0)$. By a result due to
Abhyankar, see for example \cite{Aby}, $ord_{y}(F(y,\bar 0,\bar
0))=\Sigma_{x} I(C_{\bar u^{0}},C_{\bar v^{0}}, (x0))$ at common
solutions $(x,0)$ to $C_{\bar u^{0}}$ and $C_{\bar v^{0}}$ over
$y=0$. As $C_{\bar u^{0}}$ and $C_{\bar v^{0}}$ are Weierstrass
polynomials in $x$, this is just $I(C_{\bar u^{0}},C_{\bar
v^{0}},(0,0))$. By the previous lemma and the fact that $F(y,\bar
u,\bar v)$ is generically reduced (see argument (\dag) below), it
is therefore sufficient to prove that the Zariski multiplicity of
the cover $(*)$ at $(00,\bar u^{0},\bar v^{0})$ equals the Zariski
multiplicity of the cover $Spec(L[y,\bar u',\bar v']/<F(y,\bar
u',\bar v')>)\rightarrow Spec(L[\bar u',\bar v'])$ (**) at
$(0,\bar 0,\bar 0)$. Suppose the Zariski multiplity of $(**)$
equals $n$. Then there exist distinct
$y_{1},\ldots,y_{n}\in\mathcal{V}_{0}$ and $(\bar{\delta},
\bar{\epsilon})$ generic in ${\mathcal V}_{(\bar 0,\bar 0)}\cap
U\times V$ such that $F(y_{i},\bar{\delta},\bar{\epsilon})$ holds.
Consider $Q(\bar u',\bar v')=Res(F(y,\bar u',\bar v'),{\partial
F/\partial y}(y,\bar u',\bar v'))$. This defines a closed subset
of $U\times V$ defined over $L$, we claim that this in fact proper
closed $(\dag)$. By the fact that the family is sufficently
deformable at $(\bar u_{0},\bar v_{0})$, we can find $(\bar u,\bar
v_{0})$ such that $C_{\bar u}$ intersects $C_{\bar v_{0}}$
transversally at simple points. Without loss of generality, making
a linear change of coordinates, we may suppose that for there do
not exists points of intersection of the form $(x_{1}y)$ and
$(x_{2}y)$ for $x_{1}\neq x_{2}$. By Abhyankar's result, this
implies that $F(y,\bar u',\bar v_{0})$ has no repeated roots.
Then, by genericity of $(\bar\delta,\bar\epsilon)$, we have that
$Q(\bar\delta,\bar\epsilon)\neq 0$. Hence $F(y_{i},\bar
\delta,\bar \epsilon)$ is a non-repeated root. By Abhyankar's
result, we can find a unique $x_{i}$ with $(x_{i}y_{i})$ a common
solution to the deformed curves $C_{\bar u^{0}}^{\bar\delta}$ and
$C_{\bar v^{0}}^{\bar\epsilon}$. We claim that each
$(x_{i}y_{i})\in {\mathcal V}_{00}$. As $C_{\bar
u^{0}}^{\bar\delta}(x_{i}y_{i})=0$, by the fact $(\bar
u^{0},\bar\delta,y_{i})$ specialises to $(\bar u^{0},\bar 0,0)$
and $C_{\bar u^{0}}$ is a Weierstrass polynomial in $x$, we have
that $\pi(x_{i})=0$ as well. This shows that the Zariski
multiplicity of the cover $(*)$ is at least $n$. A virtually
identical argument shows that the Zariski multiplicity of the
cover $(*)$ is at most $n$ as well.

\end{proof}
We now have the following result;\\

\begin{lemma}

Let $C_{\bar u^{0}}$ and $C_{\bar v^{0}}$ be reduced curves,
having finite intersection, then the Zariski multiplicity of the
cover $(*)$ at $((0,0),\bar u^{0},\bar v^{0})$ equals the
intersection multiplicity $I(C_{\bar u^{0}},C_{\bar v^{0}},(0,0))$
of $C_{\bar u^{0}}$ and $C_{\bar v^{0}}$ at $(0,0)$.

\end{lemma}

\begin{proof}

We have $C_{\bar u^{0}}=s(u^{0}_{ij},x,y)$ and $C_{\bar
v^{0}}=t(v^{0}_{ij},x,y)$. By making the substitutions $\bar
U=\bar u^{0}+\bar u$ and $\bar V=\bar v^{0}+\bar v$, we may assume
that $\bar u^{0}=\bar v^{0}=\bar 0$. Moreover, we can suppose
that;\\

 $s(\bar 0_{ij},x,0)\neq 0$ and $t(\bar 0_{ij},x,0)\neq 0$. (**)\\

 This can be achieved by making the invertible linear change of variables
 $(x'=x,y'=\lambda x+\mu y)$ with $(\lambda,\mu) \in L^{2}$ and $\mu\neq 0$,
  noting that as $C_{\bar u_{0}}$ and $C_{\bar v_{0}}$ are curves,
  for some choice of $(\lambda,\mu)$, the corresponding polynomials $s(u^{0}_{ij},x,y)$ and
  $t(v^{0}_{ij},x,y)$ do not vanish identically on the line $\lambda x+\mu y=0$. It is trivial to
  check that the transformation preserves both Zariski multiplicity and intersection
  multiplicity, so our calculations are not effected.\\

  We may then apply the Weierstrass preparation theorem, Lemma 3.1, in the ring
$L[[u_{ij},v_{ij},x,y]]$, obtaining factorisations
$s(u_{ij},x,y)=U_{1}(u_{ij},x,y)S(u_{ij},x,y)$ and
$t(v_{ij},x,y)=U_{2}(v_{ij},x,y)T(v_{ij},x,y)$ where $U_{1}$ and
$U_{2}$ are units in the local rings $L[[u_{ij},x,y]]$ and
$L[[v_{ij},x,y]]$, $S,T$ are Weierstrass polynomials in $x$ with
coefficients in $L[[u_{ij},y]]$ and $L[[v_{ij},y]]$ respectively.
A close inspection of the Weierstrass preparation theorem, see
\cite{Aby1}, shows that we can obtain the following uniformity in
the parameters $\bar u$ and $\bar v$;\\

Namely, if $U=\{u_{ij}:s(u_{ij},x,0)\neq 0\}$ and
$V=\{v_{ij}:t(v_{ij},x,0)\neq 0\}$, are the constructible sets for
which $(*)$ holds, then if we let $R_{U}$ and $R_{V}$ denote the
coordinate rings of $U$ and $V$, we may assume $U_{1},U_{2}$ lie
in $R_{U}[[x,y]]$ and the coefficients of $S,T$ lie in
$R_{U}[[y]]$ and $R_{V}[[y]]$ respectively. By Lemma 3.2, we may
assume that $U_{1},U_{2},S$ and $T$ lie in a finite etale
extension $R_{U\times V}[x,y]^{ext}$ of the algebra $A=R_{U\times
V}[x,y]$ (again, possibly after localisation corresponding
to an open subvariety of $Spec(A)$. Now we have the sequence of morphisms;\\

$Spec({R_{U\times V}[x,y]^{ext}\over <U_{1}S,U_{2}T>})\rightarrow
Spec({R_{U\times V}[x,y]\over <s,t>})\rightarrow Spec(R_{U\times
V})$. \\

We claim that the left hand morphism is etale at the point $(\bar
0,\bar 0,(00)^{lift})$. This follows from the fact that
$R_{U\times V}[x,y]^{ext}$ is an etale extension of $R_{U\times
V}[x,y]$ and the maximal ideal given by $(\bar 0,\bar
0,(00)^{lift})$ contains $<U_{1}S,U_{2}T>$. Now consider the
cover;\\

$Spec({R_{U\times V}[x,y]^{ext}\over <U_{1}S,U_{2}T>})\rightarrow
Spec(R_{U\times V})$ (***)\\

For $\bar u,\bar v$ in $U\times V$, the fibre of this cover over
$\bar u,\bar v$ corresponds exactly to the intersection of the
reducible curves $C'_{\bar u}$ and $C'_{\bar v}$ which lift the
original curves $C_{\bar u}$ and $C_{\bar v}$ to an etale cover of
$Spec(L[xy])$. By Theorem 1.8 and Lemma 2.3 (in the case when
$C_{\bar u_{0}},C_{\bar v_{0}}$ intersect at simple points) or
Lemma 2.7 (for singular points of intersection) and the
corresponding Lemma 2.9 for reducible covers, it is sufficent to
show that the Zariski multiplicity of the cover (***) at $(\bar
0,\bar 0,(00)^{lift})$ corresponds to the intersection
multiplicity of the curves $C'_{\bar u_{0}},C'_{\bar v_{0}}$ at
$(00)^{lift}$. The idea now is to apply Lemma 4.7 to the
Weierstrass factors of $C'_{\bar u}$ and $C'_{\bar v}$. This will
be achieved by the "unit removal" lemma below (Lemma 4.15).

\end{proof}

We first require some more definitions and a moving lemma for curves;\\

\begin{defn}

Let $X\rightarrow Spec(L[x,y])$ be an etale cover in a
neighboorhood of $(0,0)$, with distiguished point $(0,0)^{lift}$.
We call a curve $C$ on $X$ passing through $(0,0)^{lift}$
Weierstrass if, in the power series ring $L[[x,y]]$, the defining
equation of $C$ may be written as a Weierstrass polynomial in $x$
with coefficients in $L[[y]]$.

\end{defn}

\begin{defn}

Let $F\rightarrow U\times V$ be a finite equidimensional cover of
a smooth base of parameters $U\times V$ with a section $s:U\times
V\rightarrow F$. We call the cover Weierstrass with units if the
fibres $F(\bar u,\bar v)$ can be written as the intersection of
reducible curves $C_{\bar u}'$ and $C_{\bar v}'$ in an etale cover
$A_{\bar u,\bar v}$ of $U_{\bar u,\bar v}\subset Spec(L[x,y])$
with the distinguished point $s(\bar u,\bar v)$ lying above
$(0,0)$ and $C_{\bar u}', C_{\bar v}'$ factoring as $U_{\bar
u}F_{\bar u}$ and $U_{\bar v}F_{\bar v}$ with $U_{\bar u},U_{\bar
v}$ units in the local ring $O_{s(\bar u,\bar v),A_{\bar u,\bar
v}}$ and $F_{\bar u},F_{\bar v}$ Weierstrass
curves in $A_{\bar u,\bar v}$.\\

Let hypotheses on $F,U$ and $V$ be as above. We call the cover
Weierstrass if the fibres $F(\bar u,\bar v)$ can be written as
above but with $C'_{\bar u},C'_{\bar v}$ Weierstrass curves in
$A_{\bar u,\bar v}$.

We say that a Weierstrass cover (with units) factors through the
family of projective degree $d$ and degree $e$ curves if the cover
$F\rightarrow U\times V$ factors as $F\rightarrow F'\rightarrow
U\times V$ where $F'\rightarrow U\times V$ is the finite
equidimensional cover obtained by intersecting the families
$Q_{d}$ and $Q_{e}$ restricted to $U$ and $V$.\\

\end{defn}

\begin{lemma}

The cover (***) in Lemma 4.8 is a Weierstrass cover with units
factoring through the family of projective degree $d$ and degree
$e$ curves.

\end{lemma}

\begin{proof}

Clear by the above definitions.

\end{proof}

\begin{lemma}{Moving Lemma for Reduced Curves}

Let $Q_{d}$ and $Q_{e}$ be the families of \emph{all} projective
degree $d$ and degree $e$ curves. That is, with the usual
coordinate convention $x=X/Z,y=Y/Z$, $Q_{d}$ consists of all
curves of the form $s(\bar u,x,y)=\sum_{0\leq i+j\leq
d}u_{ij}x^{i}y^{j}$. Then, if $\bar u,\bar v$ are chosen in $L$,
so that the reduced curves $C_{\bar u}$ and $C_{\bar v}$ are
defined over $L$, if the tuple $\bar u'$ is chosen to be generic
in $U$ over $L$, the deformed curve $C_{\bar u}^{\bar u'}$
intersects $C_{\bar v}$ transversally at simple points.

\end{lemma}

\begin{proof}

We can give an explicit calculation;\\

Let $C_{\bar u}^{\bar u'}$ be defined by the equation $s(\bar
u',x,y)=\sum_{0\leq i+j\leq d}u_{ij}'x^{i}y^{j}$ and $C_{\bar v}$
by $t(\bar v,x,y)=\sum_{0\leq i+j\leq e}v_{ij}x^{i}y^{j}$ with
$\{v_{ij}:0\leq i+j\leq e\}\subset L$ and $\{u_{ij}':0\leq i+j\leq
d\}$ algebraically independent over $L$. Let $(x_{0}y_{0})$ be a
point of intersection, then $dim(x_{0}y_{0}/L)=1$, otherwise
$dim(x_{0}y_{0}/L)=0$ and, as $L$ is algebraically closed, we must
have that $x_{0},y_{0}\in L$. Substituting $(x_{0}y_{0})$ into the
equation $s(\bar u',x,y)=0$, we get a non trivial linear
dependence over $L$ between $u'_{00}$ and $u'_{ij}$ for $1\leq
i+j\leq d$ which is impossible. Now, the locus of singular points
for $C_{\bar v}$ is defined over $L$ and hence $(x_{0}y_{0})$ is a
simple point of $C_{\bar v}$. Now we further claim that $s(\bar
u',x,y)=0$ defines a non-singular curve in $P^{2}(K_{\omega})$
with transverse intersection to $C_{\bar v}$ Consider the
conditions Sing($\bar u$) given by $\exists x_{0}\exists
y_{0}(({\partial s\over\partial x}(x_{0}y_{0}),{\partial
s\over\partial y}(x_{0}y_{0}))=(0,0))$ and Non-Transverse($\bar
u$) by $\exists x_{0}\exists y_{0}({\partial s\over\partial
x}(x_{0}y_{0}){\partial t\over\partial y}(x_{0}y_{0})-{\partial
s\over\partial y}(x_{0}y_{0}){\partial t\over\partial
x}(x_{0}y_{0})=0)$ By the properness of $P^{2}(K_{\omega})$, these
conditions define closed subsets  of the parameter space $U$
defined over $L$. We claim that this in fact a proper closed
subset. This can be proved in a number of ways. In the case where
we restrict ourselves to affine curves, the result follows from  a
classical result of Kleiman, see \cite{Ha}, as affine space
$A^{2}(K_{\omega})$ is homogenous for the action of the additive
group $(A^{2}(K_{\omega}),+)$. More generally, we can use the
moving lemma, given in \cite{Fu}, by observing that the class of
all degree $d$ projective curves is closed under rational
equivalence. We can
also use the following enumerative method, we will abbreviate $K_{\omega}$ to $K$;\\

Consider the vector bundle $S$ of dimension $3$ on $P^{2}(K)$
given by $H^{0}(O_{P^{2}}(K)(d))/(I_{p}^{2}(d))$ where
$I_{p}^{2}(d)$ is the vector space of degree $d$ projective curves
vanishing to second order at $p$. Let $F_{1}$ be the projective
equation of $C_{\bar u}$ and $F_{2}$ a generically independent
projective curve of degree $d$. Then we have a map from the
trivial vector bundle $V_{2}$ of dimension $2$ on $P^{2}(K)$ into
$S$ given by $\Phi_{p}:\lambda F_{1}\oplus \mu F_{2}\rightarrow
S$. The number of singular curves in the pencil is given by the
$\{p\in P^{2}(K):\Phi_{p}$ has non-trivial kernel $\}$. This is
exactly the second chern class $Ch_{2}(S)$ which has codimension
$2$ as a cycle on $P^{2}(K)$. This method gives infinitely many
non-singular curves projective curves of degree $d$ in the pencil.\\

For the transversality calculation, we first assume that $C_{\bar
v}$ is non-singular. Consider the vector bundle $S$ of dimension
$2$ on $C_{\bar v}$ given by
$H^{0}(O_{P^{2}}(K)(d)))/(Tan_{p}(d))$ where $Tan_{p}(d)$ is the
vector space of degree $d$ projective curves tangent to $C_{\bar
v}$ at $p$. Again we have a map from the trivial vector bundle
$V_{2}$ of dimension $2$ on $C_{\bar v}$ given as above by
$\Phi_{p}$. The number of curves in the pencil which are tangent
to $C_{\bar v}$ is the first chern class $Ch_{1}(S)$ which has
codimension $1$ as a cycle on $C_{\bar v}$. Again this gives an
infinite number of curves with transverse intersection in the
pencil. In case $C_{\bar v}$ is singular, we obtain a jump in rank
of $S$ at the finitely many singular points, so cannot apply the
above method. We resolve this as follows, let
$\{p_{1},\ldots,p_{r}\}$ be the finitely many singular points on
$C_{\bar v}$. Choosing $2(d+1)$ generically independent
non-singular points on $C_{\bar v}$, we can find $2$ independent
degree $d$ projective curves $F_{1}$ and $F_{2}$ not passing
through the singular set of $C_{\bar v}$. Now we blow up
$P^{2}(K)$ along the curve $C_{\bar v}$. In the case when
$p_{1},\ldots,p_{r}$ are all not cusp points, the exceptional
divisor $E$ of the blow up $\pi:Bl\rightarrow P^{2}(K)$ consists
of a smooth curve $C_{0}$ and $r$ copies of $P_{1}$. We label the
crossings $\{q_{1},q_{2},\ldots,q_{2r}\}$ on $C_{0}$, mapping to
the singular points of $C_{\bar v}$. We obtain a pencil of curves
on $Bl$ by taking $(\lambda F_{1}+\mu F_{2})'$ the proper
transform of the curves $\lambda F_{1}+\mu F_{2}$ in $Bl$ (by
construction, these curves avoid $\{q_{1},\ldots,q_{2r}\}$. Now we
consider the following vector bundle $S$ on $C_{0}$. At
$C_{0}\setminus \{q_{1},\ldots q_{2r}\}$, this is just the pull
back $\pi^{*}S$ of the bundle considered above. At the crossings
$(q_{2j-1},q_{2j})$, it consists of the spaces
$H^{0}(O_{P^{2}}(K)(d))/V^{+}_{p_{j}})$ and
$H^{0}(O_{P^{2}}(K)(d))/V^{-}_{p_{j}})$, where $V^{+}_{p_{j}}$ and
$V^{-}_{p_{j}}$ consist of the spaces of projective curves of
degree $d$ passing through $p_{j}$ and tangent to the principal
axes of the normal cone at $p_{j}$. An easy calculation shows that
$S$ is a vector bundle of rank $2$ on the smooth projective curve
$C_{0}$. Now consider the trivial bundle $V^{2}$ on $C_{0}$ given
by $\lambda F_{1}+\mu F_{2}$ and the corresponding map $\Phi_{p}$.
Again, the points on $C_{0}$ for which this has non-zero kernel is
given by $Ch_{1}(S)$. This gives infinitely many curves in the
family $\lambda F_{1}+\mu F_{2}$ transverse to $C_{0}$, the
corresponding curves on $P^{2}(K)$ are then transverse to $C_{\bar
v}$ as required. In the case when $\{p_{1},\ldots,p_{r}\}$
contains cusps, we perform finitely many blow ups to obtain an
exceptional divisor of smooth curves with normal crossings. A
similar calculation (omitted), using the same method, works.

\end{proof}

\begin{rmk}

If we restrict the family of curves, the result in general fails.
A simple example is given by the family of all projective degree
$3$ curves $Q_{3}^{0,0}$ passing through $(0,0)$ with $x=X/Z$ and
$y=Y/Z$. If we take $C_{\bar v}$ to be the cusp $x^{2}-y^{3}$,
then any curve in $Q_{3}^{0,0}$ will have a non-transverse
intersection with $C_{\bar v}$ at the origin. In general we have
to use deformation theory arguments or enumerative methods to
decide this question.

\end{rmk}

\begin{lemma}{Moving Lemma for Curves with Finitely Many Marked
Points}

Let hypotheses be as in the previous lemma with $C_{\bar u}$ and
$C_{\bar v}$ defining reduced curves. Suppose also that there
exists finitely many marked points $\{p_{1},\ldots,p_{n}\}$  on
$C_{\bar v}$ defined over $L$. Then for $\bar u'\in U$ generic
over $L$ the deformed curve $ C_{\bar u}^{\bar u'}$ intersects
$C_{\bar v}$ transversely at finitely many simple points excluding
the set $\{p_{1},\ldots,p_{n}\}$.

\end{lemma}

\begin{proof}

As before, the condition that $\bar u'$ defines a curve $C_{\bar
u}^{\bar u'}$ either with non-transverse intersection to $C_{\bar
v}$ or passing through at least one of the points
$\{p_{1},\ldots,p_{n}\}$ is a closed subset of $U$ defined over
$L$. Using the above proof and the obvious fact that we can find a
curve $C_{\bar u}^{\bar u'}$ not passing through any of the points
$\{p_{1},\ldots,p_{n}\}$, we see that it is proper closed.

\end{proof}

\begin{lemma}{Unit Removal for Reduced Curves}

Let $(\pi,s):F\rightarrow U\times V$ be a Weierstrass cover with
units factoring through projective degree $d$ and degree $e$
curves. Let $(\bar u,\bar v)\in U\times V$, then there exists a
Weierstrass cover $(\pi',s'):F^{-}\rightarrow U'\times V'$ with
$U'\subset U$ and $V'\subset V$ open subsets, $(\bar u,\bar v)\in
U'\times V'$, such that $Mult_{(\bar u,\bar v,s(\bar u,\bar
v))}(F/U\times V)=Mult_{\bar u,\bar v,s'(\bar u,\bar
v))}(F^{-}/U'\times V')$.

\end{lemma}

\begin{proof}

Let $C_{\bar u}'$ and $C_{\bar v}'$ be the Weierstrass curves with
units in $A_{\bar u,\bar v}$ lifting the curves $C_{\bar u}$ and
$C_{\bar v}$. Now suppose that $Mult_{\bar u,\bar v,s(\bar u,\bar
v)}(F/U\times V)=n$. Then we can find $(\bar u',\bar
v')\in{\mathcal V}_{\bar u \bar v}\cap U\times V$ generic over $L$
such that the deformed curve $C_{\bar u}^{\bar u'}$ intersects
$C_{\bar v}^{\bar v'}$ at the $n$ distinct points
${x_{1},\ldots,x_{n}}$ in ${\mathcal V}_{s(\bar u,\bar v)}$. Now
using the Weierstrass factorisations of $C_{\bar u}^{\bar u'}$ and
$C_{\bar v}^{\bar v'}$, we claim that $U_{\bar u}^{\bar
u'}(x_{i})\neq 0$ and $U_{\bar v}^{\bar v'}(x_{i})\neq 0$. Suppose
not, then $U_{\bar u}^{\bar u'}(x_{i})=U_{\bar v}^{\bar
v'}(x_{i})=0$ and as $(\bar u',\bar v',x_{i})$
   specialises to $(\bar u,\bar v, s(\bar u,\bar v))$, then $U_{\bar
   u}(s(\bar u,\bar v))=U_{\bar v}(s(\bar u,\bar v))=0$. This
   contradicts the fact that $U_{\bar u}$ and $U_{\bar v}$ are
   units in the local ring $O_{s(\bar u,\bar v),A_{\bar u,\bar
   v}}$. Therefore, we must have that $F_{\bar u}^{\bar u'}(x_{i})=
   F_{\bar v}^{\bar v'}(x_{i})=0$. This shows that
    $Mult_{\bar u,\bar v,s(\bar u,\bar v)}(F^{-}/U\times V)\geq n$
   where $F^{-}\rightarrow U\times V$ is the cover of $U\times V$
   obtain by taking as fibres $F^{-}(\bar u,\bar v)$ the intersection
   of the Weierstrass factors $F_{\bar u}$ and $F_{\bar v}$.
   Formally, if $F$ is defined by $Spec({R_{U\times
   V}[x,y]^{ext}\over <U_{1}S,U_{2}T>})$ then $F^{-}$ is defined
   by $Spec({R_{U\times V}[x,y]^{ext}\over <S,T>})$.
   Clearly as $F^{-}\subset F$ is a union of components of $F$, we
   have that $Mult_{\bar u,\bar v,s(\bar u,\bar
   v)}(F^{-}/U\times V)\leq n$ as well. This proves the lemma.

\end{proof}

We now complete the proof of Lemma 4.8. By unit removal, it is
sufficient to compute the Zariski multiplicity of the cover\\

$Spec({R_{U\times V}[x,y]^{ext}\over <S,T>})\rightarrow
Spec(R_{U\times V})$\\

The fibre over $(\bar u,\bar v)$ of this cover corresponds exactly
to the intersection of the Weierstrass curves $F_{\bar u}$ and
$F_{\bar v}$ lifting $C_{\bar u}$ and $C_{\bar v}$. We then use
Lemma 2.7, noting that the Weierstrass factors are still reduced,
see \cite{Aby1}, to finish the result, with the straightforward
modification that we work in a uniform family of etale covers.\\

We now turn to the problem of non-reduced curves. We will show the
following stronger version of Lemma 4.8\\

\begin{lemma}

Let $C_{\bar u^{0}}$ and $C_{\bar v^{0}}$ be non-reduced curves
having finite intersection, then the Zariski multiplicity of the
cover (*) at $((0,0),\bar u^{0},\bar v^{0})$ equals the
intersection multiplicity $I(C_{\bar u^{0}},C_{\bar v^{0}},(0,0))$
of $C_{\bar u^{0}}$ and $C_{\bar v^{0}}$ at $(0,0)$.

\end{lemma}

First, we will require some more lemmas.

\begin{lemma}

Let $C_{\bar u_{0}}$ and $C_{\bar v_{0}}$ be reduced curves
intersecting transversally at $(0,0)$. Then the Zariski
multiplicity, left multiplicity and right multiplicity of the
cover (*) at $((0,0),\bar u^{0},\bar v^{0})$ equals $1$.

\end{lemma}

\begin{proof}

First note that by Lemma 2.6 (and corresponding Lemma 2.9), and
the fact that a generic deformation $C_{\bar v_{0}}^{\bar v'}$
will still intersect $C_{\bar u_{0}}$ transversally by Lemma 4.12,
it is sufficient to prove the result for right multiplicity.\\

In order to show this we require the following result, given for
analytic curves in \cite{Aby1}, we will only need the result for polynomials;\\

Implicit Function Theorem:\\

       If $G(X,Y)$ is a power series with
$G(0,0)=0$ then $G_{Y}(0,0)\neq 0$ implies there exists a power
series $\eta(X)$ with $\eta(0)=0$ such that $G(X,\eta(X))=0$.\\

In order to show that $Right Mult_{(0,0),\bar u^{0},\bar
v^{0}}(F'/U\times V)=1$, where $F'$ is the family obtained by
intersecting degree $d$ and degree $e$ curves, we apply the
implicit function theorem to the curve $C_{\bar u^{0}}$ at the
point $(0,0)$ of intersection with $C_{\bar v^{0}}$. Let $G(X,Y)$
and $H(X,Y)$ denote the polynomials defining the curves. We have
that $G(0,0)=H(0,0)=0$. Moreover, as the first curve is
non-singular at $(0,0)$, we may also assume that $G_{Y}(0,0)\neq
0$. Now let $\eta(X)$ be given by the theorem. As the intersection
of the curves $C_{\bar u^{0}}$ and $C_{\bar v^{0}}$ is transverse,
$ord_{X}H(X,\eta (X))=1$.Now we
have the sequence of maps;\\

$L[\bar v]\rightarrow {L[X,Y][\bar v]\over <G(u^{0},X,Y),H(\bar
v,X,Y)>}\rightarrow {L[X]^{ext}[Y][\bar v]\over <Y-\eta(X),H(\bar
v,X,Y)>}$.\\

where $L[X]^{ext}$ is an etale extension of $L[X]$ containing
$\eta(X)$. (Note that $\eta(X)$ is trivially algebraic over
$L(X)$). This corresponds to a sequence of finite covers
$F_{1}\rightarrow F'(u_{0},V)\rightarrow Spec(L[\bar v])$. The
left hand morphism is trivially etale at $(\bar
v^{0},(00)^{lift})$, hence it is sufficient to compute the Zariski
multiplicity of $F'\rightarrow Spec(L[\bar v])$ at $(\bar v^{0},
(00)^{lift})$ by Lemma 2.3 (or corresponding Lemma 2.9). This is a
straightforward calculation, the fibre over $\bar v^{0}$ consists
of the scheme $Spec({L[X,\eta(X)]\over G(X,\eta(X))})=Spec(L)$ as
$ord_{X}(H(X,\eta(X)))=1$, hence is etale at the point $(\bar
v^{0},(00)^{lift})$. By Theorem 1.4, the Zariski multiplicity is $1$.\\

\end{proof}

\begin{lemma}

Let hypotheses be as in Lemma 4.17, then for any $(\bar u',\bar
v')\in {\mathcal V_{(\bar u^{0},\bar v^{0})}}$, we have that $Card
(F'(\bar u',\bar v')\cap {\mathcal V_{0,0}})=1$

\end{lemma}

\begin{proof}

This follows immediately from Lemma 4.17 and Lemma 2.4.

\end{proof}

\begin{defn}

For ease of notation, given curves $C_{\bar u}$ and $C_{\bar v}$
of degree $d$ and degree $e$ intersecting at $x\in
P^{2}(K_{\omega})$, we define $Mult_{x}(C_{\bar u},C_{\bar v})$ to
be the corresponding Zariski multiplicity of the cover
$F'\rightarrow U\times V$ at the point $(x,\bar u,\bar v)$.
Similarly for left/right multiplicity.

\end{defn}

We can now give the proof of Lemma 4.16;\\

\begin{proof}

Case 1. $C_{\bar v_{0}}$ is a reduced curve (possibly having
components). Write $C_{\bar u^{0}}$ as $G_{1}^{n_{1}}(X,Y)\ldots
G_{m}^{n_{m}}(X,Y)=0$ with $G_{i}$ the reduced irreducible
components of $C_{\bar u_{0}}$ with degree $d_{i}$ passing through
$(0,0)$. Choose $\bar \epsilon_{1}^{1},\ldots \bar
\epsilon_{1}^{n_{1}},\ldots \bar \epsilon_{i}^{j},\ldots \bar
\epsilon_{m}^{1},\ldots, \bar \epsilon_{m}^{n_{m}}$ independent
generic in $U_{i}$, the parameter space for degree $d_{i}$
projective curves with $\bar \epsilon_{i}^{j}\in {\mathcal V_{\bar
u_{i}^{0}}}$, where $\bar u_{i}^{0}$ defines $G_{i}$. By repeated
application of Lemma 4.14, the deformed curves $G_{i}^{\bar
\epsilon_{i}^{j}}=0$ intersect $C_{\bar v_{0}}$ transversely at
disjoint sets of points We denote by $Z_{\bar \epsilon_{i}^{j}}$
those points lying in ${\mathcal V_{00}}$.  Now the curve defined
by $\prod_{ij}G_{i}^{\bar \epsilon_{i}^{j}}=0$ is a deformation
$C_{\bar u^{0}}^{\bar \epsilon}$ of $C_{\bar u^{0}}$. We let
$Z_{\bar \epsilon}$ denote the points of intersection of $C_{\bar
u^{0}}^{\bar \epsilon}$ with $C_{\bar v_{0}}$ in ${\mathcal V}_{00}$. Then we have;\\

$Z_{\bar \epsilon}=\bigcup_{ij} Z_{\bar \epsilon_{i}^{j}}$\\

$Card(Z_{\bar \epsilon})=\sum_{ij}Card(Z_{\bar \epsilon_{i}^{j}})$\\

By Lemma 2.4, we have that\\

$Left Mult_{(00)}(C_{\bar u^{0}},C_{\bar v^{0}})=\sum_{x\in
Z_{\bar\epsilon}}Left Mult_{x}(C_{\bar u^{0}}^{\bar
\epsilon},C_{\bar v^{0}})$\\
$=\sum_{i,j}\sum_{x\in Z_{\bar\epsilon_{i}^{j}}}Left
Mult_{x}(C_{\bar
u^{0}}^{\bar\epsilon},C_{\bar v^{0}})$ (*)\\

We now claim that for a point $x\in Z_{\bar \epsilon_{i}^{j}}$,\\

$Left Mult_{x}(C_{\bar u^{0}}^{\bar \epsilon},C_{\bar v^{0}})=Left
Mult_{x}(G_{i}^{\bar \epsilon_{i}^{j}},C_{\bar v_{0}})$ (**)\\

 This follows as both the reduced curves $C_{\bar
u_{0}}^{\bar \epsilon}$ and $G_{i}^{\bar \epsilon_{i}^{j}}$
intersect $C_{\bar v_{0}}$ transversely at $x$. Hence, in both
cases the left multiplicity is $1$, by Lemma 4.17.\\

Combining $(*)$ and $(**)$, we obtain;\\

$Left Mult_{(00)}(C_{\bar u^{0}},C_{\bar
v^{0}})=\sum_{i,j}\sum_{x\in Z_{\bar \epsilon_{i}^{j}}}Left
Mult_{x}(G_{i}^{\bar \epsilon_{i}^{j}},C_{\bar v^{0}})$\\

Now using Lemma 2.4 again gives that;\\

$Left Mult_{(00)}(C_{\bar u^{0}},C_{\bar
v^{0}})=\sum_{i=1}^{m}n_{i}Left Mult_{(00)}(G_{i},C_{\bar v^{0}})$ (***)\\

 If we go through exactly the same calculation with Mult replacing
Left Mult, we see as well that\\

$Mult_{(00)}(C_{\bar u^{0}},C_{\bar
v^{0}})=\sum_{i=1}^{m}n_{i}Mult_{(00)}(G_{i},C_{\bar v^{0}})$\\

By Lemma 4.8, this gives \\

$Mult_{(00)}(C_{\bar u^{0}},C_{\bar
v^{0}})=\sum_{i=1}^{m}n_{i}I(G_{i},C_{\bar v^{0}},(00))$\\

 By a straightforward algebraic calculation, see references below
for the definitive result,
this gives\\

$Mult_{(00)}(C_{\bar u^{0}},C_{\bar v^{0}})=I(C_{\bar
u^{0}},C_{\bar v^{0}},(00))$\\

as required.\\

 Case 2. Both $C_{\bar u_{0}}$ and $C_{\bar v_{0}}$ define
non-reduced curves. Write $C_{\bar u_{0}}$ as above and $C_{\bar
v_{0}}$ as $H_{1}^{e_{1}}\ldots H_{n}^{e_{n}}$ with $H_{i}$ the
reduced compoments with degree $c_{i}$ of $C_{\bar v_{0}}$ passing
through $(00)$. Then $H_{1}\ldots H_{n}=0$ defines a reduced curve
passing through $(00)$. Now repeat the argument in Case 1 for the
curves $C_{\bar u_{0}}$ and $H_{1}\ldots H_{n}=0$. Again let
$Z_{\bar \epsilon}$ be the intersection points of the deformed
curve $C_{\bar u_{0}}^{\epsilon}$ with $H_{1}\ldots H_{n}=0$ in
${\mathcal
V_{(00)}}$. By (***) of Case 1, Lemma 2.4 and Lemma 4.18 with the fact that the intersection
of $C_{\bar u^{0}}^{\bar \epsilon}$ with $H_{1}\ldots H_{n}$ is transverse, we have;\\

$Card(Z_{\bar
\epsilon})=\sum_{i=1}^{m}n_{i}Mult_{(00)}(G_{i},H_{1}\ldots
H_{n})$\\

 Now using the argument in Case 1 applied to the reduced curves
$G_{i}$ and $H_{1}\ldots H_{n}$, we have;\\

$Card(Z_{\bar
\epsilon})=\sum_{i=1}^{m}n_{i}\sum_{j=1}^{n}I(G_{i},H_{j},(00))$
(*)\\

We claim that for any component $H_{j}$\\

$Card(H_{j}\cap Z_{\bar
\epsilon})=\sum_{i=1}^{m}n_{i}I(G_{i},H_{j},(00))$\\

 This follows as the deformed curve $C_{\bar u_{0}}^{\bar \epsilon}$
 a fortiori intersects $H_{j}$ transversely at simple points. Therefore,
  again by Case 1, gives the expected multiplicity. Now, using this together
  with (*), we write $Z_{\bar \epsilon}$ as $\cup_{j}Z_{\bar \epsilon}^{j}$ where
$Z_{\bar \epsilon}^{j}$ are the disjoint sets consisting of the
intersection of $C_{\bar u_{0}}^{\bar \epsilon}$ with $H_{j}$.
Then by Lemma 2.6, we have that\\

$Mult_{(00)}(C_{\bar u^{0}},C_{\bar v^{0}})=\sum_{j}\sum_{x\in
Z_{\bar \epsilon}^{j}}Right Mult_{x}(C_{\bar u^{0}}^{\bar
\epsilon},C_{\bar v^{0}})$\\

 We can now calculate the Right Mult term by applying Case 1 to
the intersection of $C_{\bar v_{0}}$ with the reduced curve
$C_{\bar u_{0}}^{\bar \epsilon}$ at the points of intersection
$x\in Z_{\bar \epsilon}^{j}$. At a point $x\in Z_{\bar
\epsilon}^{j}$,
we have that\\

$Right Mult_{x}(C_{\bar u^{0}}^{\bar \epsilon}, C_{\bar
v^{0}})=e_{j}I(C_{\bar u^{0}}^{\bar \epsilon},H_{j},x)=e_{j}$\\

 as the intersection is transverse. Finally this gives;\\

 $Mult_{(00)}(C_{\bar u^{0}},C_{\bar
 v^{0}})=\sum_{i=1}^{m}\sum_{j=1}^{n}n_{i}e_{j}I(G_{i},H_{j},(00))$\\

 By an algebraic result, see \cite{Kir} for the case of complex
algebraic curves, or \cite{Ful}
 for its generalisation to algebraic curves in arbitrary characteristics, we have\\

 $Mult_{(00)}(C_{\bar u^{0}},C_{\bar v^{0}})=I(C_{\bar u^{0}},C_{\bar v^{0}},(00))$\\

as required.\\

\end{proof}

The following version of Bezout's theorem in all characteristics
is now an easy generalisation from the above lemma. For curves
$C_{1}$ and $C_{2}$ in $P^{2}(L)$, we let $M(C_{1},C_{2},x)$
denote the intersection multiplicity or the Zariski multiplicity,
we know from the above that the two are equivalent.

\begin{theorem}{Non-Standard Bezout}

 \ \ \ \ \ \ \ \ \ \ \ \\

Let $C_{1}$ and $C_{2}$ be projective curves of degree $d$ and
degree $e$ in $P^{2}(L)$, possibly with non-reduced components,
intersecting at finitely many points $\{x_{1},\ldots,x_{i},\ldots
x_{n}\}$,
 then we have;\\

$\sum_{i=1}^{n}M(C_{1},C_{2},x_{i})=de$.

\end{theorem}

Of course, we could just quote the algebraic result given in
\cite{Ha} (though this in fact only holds for reduced curves).
Instead we can give a non-standard proof, which in many ways is
conceptually simpler and doesn't involve any algebra;\\

\begin{proof}

Let $Q_{d}$ and $Q_{e}$ be the families of all projective degree
$d$ and degree $e$ curves. Then we have the cover $F\rightarrow
U\times V$ with $F\subset U\times V\times P^{2}(L)$ obtained by
intersecting the families $Q_{d}$ and $Q_{e}$. We have that\\

$\sum_{i=1}^{n}M(C_{1},C_{2},x_{i})=\sum_{i=1}^{n}Mult_{x_{i}\in
F(\bar u_{0},\bar v_{0})}(F/U\times V)$ \\

 where $(\bar u_{0},\bar v_{0})$ define $C_{1}$ and $C_{2}$. By
Lemma 4.3 in \cite{deP},
this equals\\

 $\sum_{x\in F(\bar u,\bar v)}Mult_{x,\bar u,\bar
v}(F/U\times V)$\\

 where $(\bar u,\bar v)$ is generic in $U\times
V$. Using, for example, the proof of Lemma 4.12, generically
independent curves $C_{\bar u}$ and $C_{\bar v}$ intersect
transversely at a finite number of simple points. Hence, by Lemma
4.17, the Zariski multiplicity calculated at these points is $1$.
As the cover $F$ has degree $de$, there is a total number $de$ of
these
points as required.\\

\end{proof}

\end{section}

\end{document}